\documentclass[aps,prl,twocolumn]{revtex4-1}

\usepackage{times}
\usepackage{mathptmx}
\usepackage{amsmath}
\usepackage{amsfonts}
\usepackage{amssymb}
\usepackage{mathrsfs}
\usepackage{enumerate}
\usepackage{bm}
\usepackage{graphicx}
\usepackage{hyperref}
\usepackage{color}
\usepackage{array}
\usepackage{booktabs}
\usepackage{verbatim}
\usepackage{pagecolor}
\usepackage{newtxtext,newtxmath}
\usepackage[greek,english]{babel}

\newcolumntype{C}{>{$}c<{$}}
\AtBeginDocument{
\heavyrulewidth=.08em
\lightrulewidth=.05em
\cmidrulewidth=.03em
\belowrulesep=.65ex
\belowbottomsep=0pt
\aboverulesep=.4ex
\abovetopsep=0pt
\cmidrulesep=\doublerulesep
\cmidrulekern=.5em
\defaultaddspace=.5em
}

\definecolor{CiteColor}{rgb}{0.4,0.58,0.79}
\definecolor{URLColor}{rgb}{0.4,0.58,0.79}
\definecolor{LinkColor}{rgb}{0.4,0.58,0.79}
\hypersetup{colorlinks,citecolor=CiteColor,urlcolor=URLColor,linkcolor=LinkColor}

\newcommand{\beq}{\begin{equation}}
\newcommand{\eeq}{\end{equation}}

\newcommand{\ue}{\mathrm{e}}
\newcommand{\upi}{\textrm{\greektext p}}
\newcommand{\ui}{\mathrm{i}}

\newcommand{\CC}{\mathbb{C}}
\newcommand{\NN}{\mathbb{N}}

\newcommand{\psqrt}{\sqrt{\phantom{\ui}}}

\newcommand{\red}{\textcolor{red}}
\newcommand{\cyan}{\textcolor{cyan}}
\newcommand{\orange}{\textcolor{orange}}
\newcommand{\green}{\textcolor{violet}}

\begin{document}

\title{Abel--Ruffini's Theorem: Complex but Not Complicated ! \texorpdfstring{\\}{}
        \normalfont{ \normalsize{\textit{A proof, using loops and roots, of the unsolvability of the quintic} }} }

\author{Paul Ramond}\email{paul.ramond@obspm.fr}
\affiliation{Laboratoire Univers et Th{\'e}ories, Observatoire de Paris, CNRS, Universit{\'e} PSL, Universit{\'e} de Paris, 92190 Meudon, France}

\begin{abstract}
\hspace{-0.17cm}

In this article, using only elementary knowledge of complex numbers, we sketch a proof of the celebrated Abel--Ruffini theorem, which states that the general solution to an algebraic equation of degree five or more cannot be written using radicals, that is, using its coefficients and arithmetic operations $+,-,\times,\div,$ and $\psqrt$. The present article is written purposely with concise and pedagogical terms and dedicated to students and researchers not familiar with Galois theory, or even group theory in general, which are the usual tools used to prove this remarkable theorem. In particular, the proof is self-contained and gives some insight as to why formulae exist for equations of degree four or less (and how they are constructed), and why they do not for degree five or more.

\end{abstract}


\maketitle

\section{Introduction}

\textit{Historical background.---}Finding a general expression for the solutions of an algebraic equation has been one of the oldest and most fruitful problems in mathematics. The history behind what was once called the ``theory of equations'' \cite{Stillwell}, is almost as rich and old as the history of mathematics itself. For example, methods for solving linear and quadratic equations have been known for at least four millennia \cite{Sesiano}, in independent places in the world. The quadratic formula taught today in school, with modern notation, was first written down by R.~Descartes in 1637 \cite{Descartes}. The introduction of the definitive $\psqrt{}$ notation (with the horizontal overbar called the ``vinculum'') was only introduced in 1525 \cite{Mazur}. Regarding cubic and quartic equations, they too had to wait until the sixteenth century to be finally solved. By then, a group of rivaling Italian mathematicians, including S.~del Ferro, N.~Tartaglia, G.~Cardano, and L.~Ferrari, made the serendipitous discovery of complex numbers while solving the general cubic equation. In 1545, a few years before their quarrels settled in a public mathematical contest \cite{Stillwell}, L.~Ferrari solved the quartic equation by reducing it to a cubic one. The quintic equation, however, would still keep these mathematicians (and all others) in check, while the idea of it being unsolvable slowly started to emerge.\\

\textit{Unsolvable equations.---}The idea of examining permutations of the solutions to study the (un-)solvability of algebraic equations dates back to the pioneering works of J.-L.~Lagrange in 1771 \cite{Lagrange}. Lagrange's ideas matured, and were finally extended to the quintic equation by P.~Ruffini in early 1800 \cite{Ruffini}. For twenty years, Ruffini tried to convince the mathematical community of the importance of his results, without success. It is only in 1821, with the help of L.-A.~Cauchy, that Ruffini's work was recognized as a stepping stone in the theory of algebraic equations. Although it turned out that Ruffini did not prove the theorem that now bears his name \textit{per se}, his results were strong enough to place serious doubt about the possibility of finding a solution to the general quintic equation. The wait was finally over in 1824 when N.-H.~Abel wrote the first complete proof of the theorem (a short proof published in 1824 \cite{Abel.small} at his own expense, and a longer, more detailed version two years later \cite{Abel.big}). His work still remained unworthy of interest to the eyes of most mathematicians, including Gauss and Cauchy themselves. Abel died aged 26 in 1829, just before his work on the unsolvability of the quintic finally received all the appreciation it deserved. He received posthumously the \textit{Grand Prix de l'Acad\'{e}mie des Sciences de Paris} in 1830, in recognition of his work. The same year also marks the publication of \'{E}.~Galois's first paper on these topics \cite{Galois}, in which he gives the premises of (now) Galois theory, a novel and elegant extension of all previous results. He too died young (aged 20 in 1832) and his work also took several decades to be fully published and recognized as revolutionary.\\

This short historical account lacks many interesting stories about these mathematicians, such as conflict of interests, encrypted communications, fatal duels, long lost and recovered memoirs, etc. The interested reader could start with J.~Sesiano's \cite{Sesiano} and J.~Stillwell's \cite{Stillwell} books and references therein for well-written and thorough presentations of these fascinating pieces of history.\\

\textit{Aim and content.---}The aim of this article is to sketch an accessible and self-contained proof of the \textit{Abel--Ruffini theorem}:
\begin{center}
\textit{No formula exists for the solution to the general} \\ \textit{equation of degree five or more, using only} \\ \textit{the operations $+,-,\times,\div,$ and $\psqrt$.}
\end{center}

The word \textit{general} is important: it emphasizes that a formula that \textit{holds for any coefficients} cannot be found. However, the theorem does not prevent \textit{some} equations to have a solution that can be written in terms of $+,-,\times,\div,\psqrt$. In most textbooks, the proof of this remarkable theorem relies on a powerful subbranch of mathematics called \textit{Galois theory}, developed quasi-exclusively by the French mathematician \'{E}.~Galois at the beginning of the nineteenth century. Galois theory solved all ``unsolvability problems'' once and for all, as well as other millennia-long problems \cite{Stillwell}. However, it is also rather advanced, usually taught in the second/third years of specialized, university-level mathematics. The first complete proof (by Abel) of the Abel--Ruffini theorem is a few years older than the birth of Galois theory. Moreover, the works of Galois took several decades to be broadly known to other mathematicians. In other words, neither Ruffini nor Abel used the methods developed by Galois to prove that some equations were unsolvable. Because it usually relies on advanced mathematics, few people in the scientific community are aware of this theorem and its underlying principles. But because Abel did not prove it this way, there must be another, perhaps simpler, way of understanding the reason why the general quintic equation does not have a solution in terms of radicals. In particular, Galois's, Abel's, and Ruffini's ideas all rely on a unique, fundamental, common point: \textit{the symmetry of an algebraic equation under the permutation of its solutions}. Based solely on this fundamental symmetry, we propose to sketch a proof of Abel--Ruffini's theorem using only elementary knowledge about complex numbers. Familiarity with complex numbers and a pen (and paper!) to draw appropriate figures are the only prerequisites to get a grasp of how the proof works. Everything else is elementary mathematics and useful notations that help present the ideas more clearly.\\

\textit{Motivation.---}The proof given here cannot be considered new. It is the result of several adaptations and simplifications of ideas that we feel compelled to attribute to the theoretical physicist B.~Katz. His ideas are presented concisely in an online video \cite{Katz}, which can be used as complementary material with dynamic illustrations. Katz's inspiration for making this video comes from a series of lectures given by physicist and mathematician V.~Arnold, which were nicely crystallized in a problems-and-solutions book published by V.~B.~Alekseev, who was Arnold's student at the time of these lectures. This book, although very well written and complete, is, however, not elementary in any sense. While Katz's video does a very good job at explaining the general idea of the proof, we found that some gaps could be filled, and some arguments could be made much simpler, especially when we get to the end of the proof. Other references dealing with the present ideas are rather scarce in the literature (academic or not). A nonexhaustive selection is located at the conclusion of the article, and can serve as complementary material to deepen one's understanding of the proof.\\

\textit{Outline.---}The remainder of this article is organized as follows. After some prerequisites and reminders regarding complex numbers are introduced, we spend some time on the quadratic equation, explaining why a quadratic formula cannot be built out of only the four basic arithmetic operations (our first impossibility result). Similar ideas are then extended successively to the cubic and quartic equation, giving stronger impossibility results at each step. By the time we get to the quintic equation, the reader should be comfortable enough with the strategy (hopefully) to see how the quadratic, cubic, and quartic cases foreshadow the proof of the the Abel--Ruffini theorem. Along the way, we also derive the cubic and quartic formulae, scarcely presented in the nonspecialized literature. Although the derivation of these formulae is interesting enough to justify their presence, they will be especially useful in light of our temporary results, and will naturally guide us step by step to Abel--Ruffini's theorem. Finally, we note that animated versions of Figures \ref{fig:examplequintic}, \ref{fig:vieta}, \ref{fig:permuquad} and \ref{fig:permucubic} are available as supplementary material for a better understanding.

\section{Prerequisites}
In this article, we are dealing with algebraic equations of degree $n\geq 2$. These equations are always of the form
\beq \label{geneq}
z^n + c_{n-1}z^{n-1} + \cdots + c_1 z + c_0 = 0 \, ,
\eeq
where $z\in\CC$ is the \textit{unknown} and the $n$ complex numbers $(c_0,...,c_{n-1})$ are the \textit{coefficients}. It is a remarkable fact, often cited as the fundamental theorem of algebra, that equation \eqref{geneq} always has exactly $n$ complex solutions. (We use \textit{solutions}, instead of \textit{roots} of polynomials, to avoid confusion with the ``root'' operation $\psqrt$ later on.) These solutions will always be denoted $(s_1,\ldots,s_n)$, and we use $s$ as a placeholder for \textit{any} of the solutions. \\

\begin{figure}[!htbp]
  \centering
	 \includegraphics[width=0.7\linewidth]{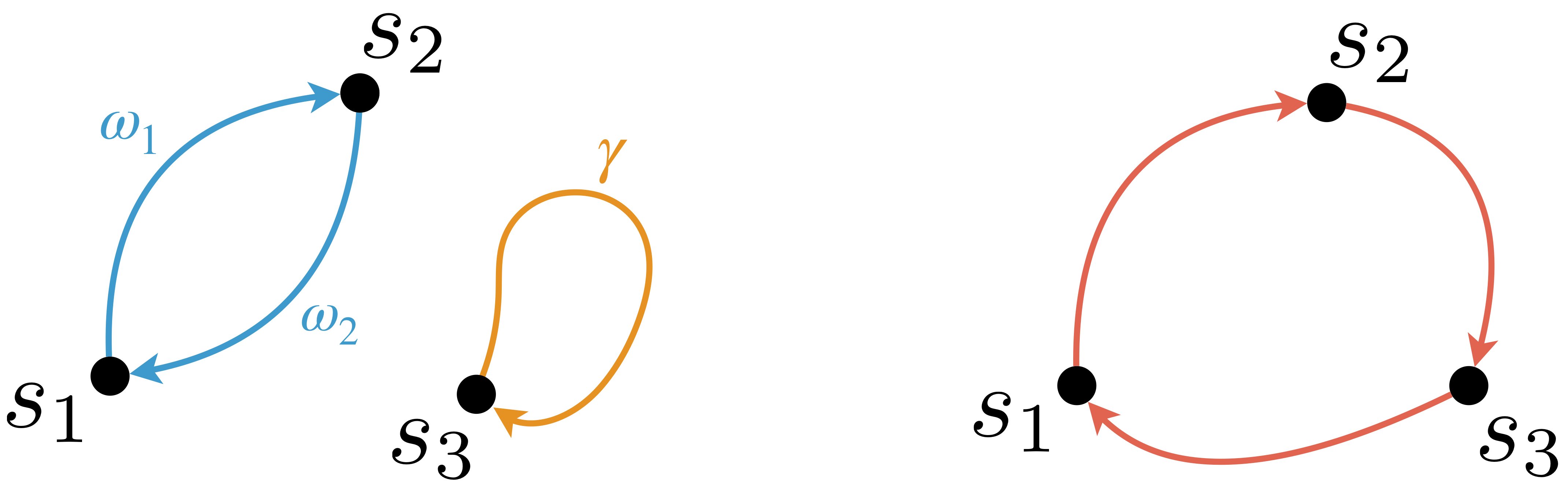}
   \caption{The paths-induced \textit{transposition} $\cyan{(}12\cyan{)}$ and \textit{cycle} $\red{(}123\red{)}$ on the solutions $(s_1,s_2,s_3)$ of some algebraic equation of degree $n\geq3$. \textit{See supplementary material for animated version.} \label{fig:permu}}
\end{figure}

\textit{Permutations.---}Our strategy will be based on picturing the solutions $(s_1,\ldots,s_n)$ in the complex plane and make them move around so as to exchange their positions, i.e., \textit{permute} them. We will need two kinds of permutation:
\begin{itemize}
    \item \textit{transpositions}, denoted $(ij)$, exchanging the position of \textit{two} solutions, i.e., $s_i\leftrightarrow s_j$. The transposition $(12)$ is depicted on the left in Fig~\ref{fig:permu},
    \item \textit{cycles}, denoted $(ijk)$, exchanging the position of \textit{three} solutions cyclically, i.e., $s_i\rightarrow s_j$, $s_j\rightarrow s_k$, and $s_k\rightarrow s_i$. The cycle $(123)$ is depicted on the right in Fig~\ref{fig:permu}.
\end{itemize}

Two permutations next to each other are to be performed successively, from \textit{left to right}. For example, $(12)(23)$ consists in exchanging $s_1$ and $s_2$, \textit{then} $s_2$ with $s_3$. Notice that the result is equivalent to the cycle $(132)$, hence there is no unique way of writing permutations. However, permutations \textit{do not commute} in general. Indeed, $(12)(23)=(132)$ and $(23)(12)=(123)$; therefore $(12)(23)\neq (23)(12)$.\\

\textit{Loops.---}One way of visualizing permutations of $(s_1,\ldots,s_n)$ is to locate them in the complex plane and make them travel along some \textit{paths}. Paths in the complex plane are just continuous curves than connect two points (we assume that they do not self-intersect, otherwise things get unnecessarily complicated). A path that \textit{closes}, i.e., connects a point to itself, is called a \textit{loop} and denoted $\gamma$, whereas a path that connects two \textit{distinct} points is simply called an \textit{unclosed path}, denoted $\omega$. These paths will be represented by arrows in our figures, and will be used to induce permutations on $(s_1,\ldots,s_n)$. For example, in Figure \ref{fig:permu} are depicted the transposition $\cyan{(}12\cyan{)}$ on the left, and the cycle $\red{(}123\red{)}$ on the right. Notice that to induce $\cyan{(}12\cyan{)}$, $s_3$ follows a loop $\orange{\gamma}$ so that only $s_1$ and $s_2$ swap places by following the unclosed paths $\cyan{\omega_1},\cyan{\omega_2}$. When speaking of permutations of solutions, we will \textit{always} imagine them traveling on these paths. \\

\begin{figure}[!htbp]
  \centering
	 \includegraphics[width=0.8\linewidth]{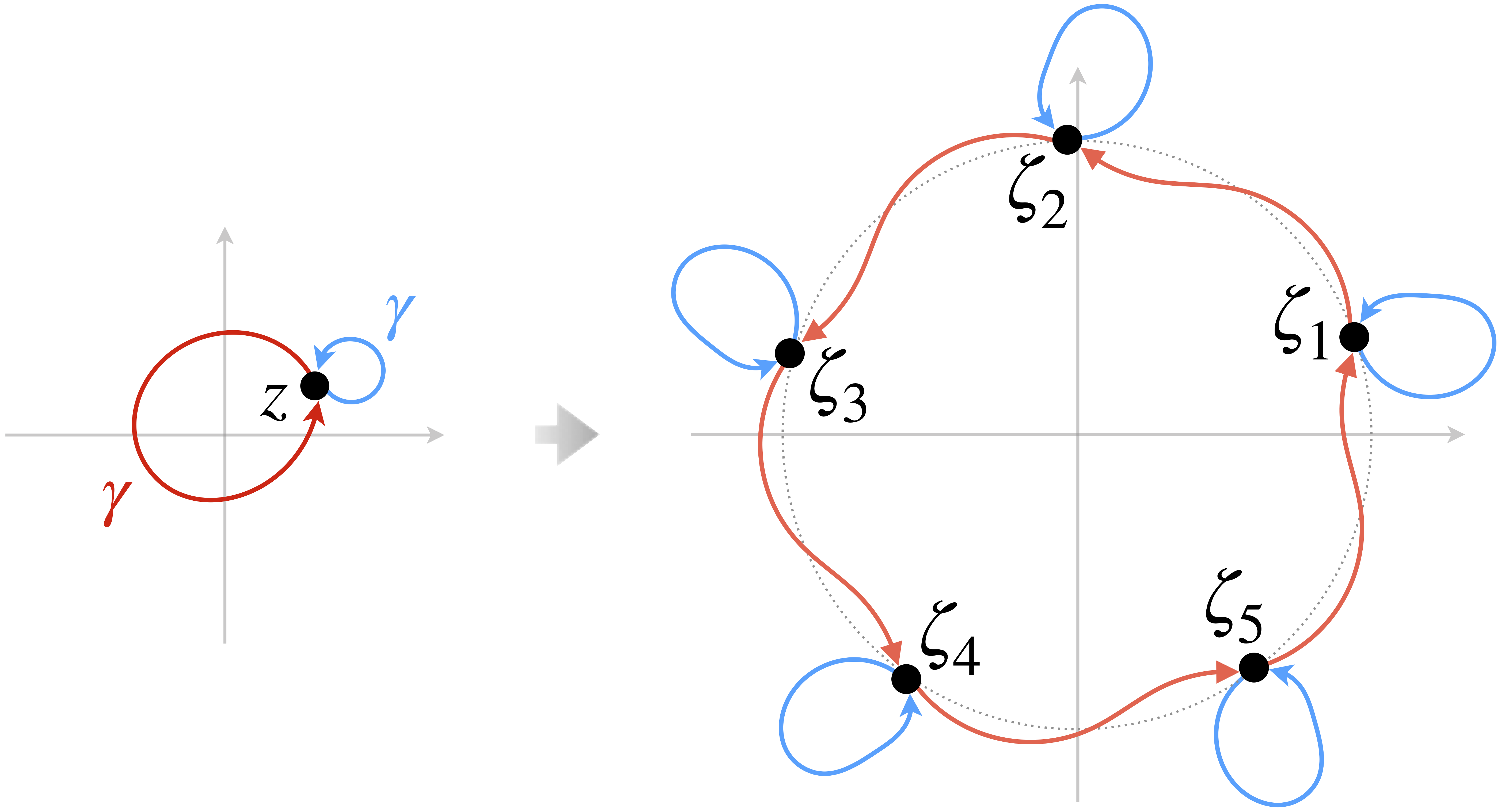}
   \caption{When $z$ follows a loop $\cyan{\gamma}$ that does not wind around the origin, its roots $\zeta_\ui$ follow loops as well (right). However, for the loop $\red{\gamma}$ that winds once, the roots then follow the red, unclosed paths. \textit{See supplementary material for animated version.} \label{fig:examplequintic}}
\end{figure}

\textit{Roots.---}Now let us examine how roots of complex numbers move around the complex plane. Fixing some complex number $z$, a \textit{root} of $z$ is some number $\zeta\in\CC$ such that $\zeta^k=z$ for some $k\in\NN$. Such a $\zeta$ is then called a $k$th root of $z$; and $z$ admits exactly $k$ such $k$th roots (this follows from the fundamental theorem of algebra).
We will deliberately use the ambiguous notation $\sqrt[k]{\phantom{\ui}}$ as a multivariable notation, i.e., for a given $k$, $\sqrt[k]{z}$ means \textit{any} of the $k$th roots of $z$. Fixing $k\in\NN$ and assuming that $z$ itself follows a loop $\gamma$, let us examine what kind of path $\sqrt[k]{z}$ follows. To this end, we use the exponential form of $z$, i.e, $z=r \ue^{\ui\theta}$ with $r=|z|$ and $\theta=\arg z$, from which we find that all $k$th roots $(\zeta_1,\ldots,\zeta_k)$ can be written explicitly as
\beq \label{roots}
\zeta_\ell=r^{1/k} \ue^{\ui(\theta +2\ell \upi )/k} \, , \quad \ell\in\{1,\ldots,k\} \, .
\eeq
From equation \eqref{roots}, one can already tell that all $k$th roots of $z$ have the same modulus : $r^{1/k}$. Geometrically, this means that they lie on the same circle (of radius $r^{1/k}$) in the complex plane. Moreover, we readily see from equation \eqref{roots} that
%
%
\beq \label{move}
\arg \zeta_\ell=\frac{\theta}{k}+ \ell \frac{2\upi}{k} \, ,
\eeq
which means that all roots are equally spaced on this circle, at angle $2\upi/k$ apart. Now suppose that $z$ goes on a journey exploring the complex plane, by traveling on a loop $\red{\gamma}$ winding once (say) around the origin, in the counterclockwise direction (in red, on the left in Figure \ref{fig:examplequintic}). As $z$ travels along $\red{\gamma}$, its $k$th roots also move around, and their position can be tracked from equation \eqref{roots} (see the red paths on the right in Figure \ref{fig:examplequintic}). Since $\red{\gamma}$ is a loop, the radius $r$ comes back to its original (i.e., pre-loop) value, and so does $r^{1/k}$. In other words: the roots remain on their circle after the path $\red{\gamma}$ (see the grey, dashed-circle on the right of Figure \ref{fig:examplequintic}). However, $\arg z$ went from $\theta$ to $\theta+2\upi$ (one counterclockwise turn). Therefore, from equation \eqref{move}, each $k$th root $\zeta_\ell$ has moved to its closest, counterclockwise neighbor, $\zeta_{\ell+1}$. In particular, the roots have followed an unclosed path. Had $z$ not wound around the origin (in blue, on the left in Figure \ref{fig:examplequintic}), its argument $\theta$ would have seen no net change after the loop $\cyan{\gamma}$, and the roots would have followed their own loops (in blue on the right in Figure \ref{fig:examplequintic}).
\\

We have seen two example of loops followed by $z$ and the result is not the same for its roots $\sqrt[k]{z}$ : in one case the roots follow a loop (blue part of Figure \ref{fig:examplequintic}), in the other they \textit{do not} (red part of Figure \ref{fig:examplequintic}). Consequently, we conclude that \textit{when $z$ follows a loop, $\sqrt[k]{z}$ does not always follow a loop}. This conclusion holds for any type of root (i.e., any $k$ in $\sqrt[k]{z}$). Since we will not need to differentiate between all these roots, we will denote by $\sqrt{z}$ \textit{any} root of $z$ (that is, \textit{any} $k$th root, whatever the value $k\in\NN$). With this notation, the takeaway result of this paragraph is simply:
\begin{center}
    \textit{When $z$ follows a loop, $\sqrt{z}$ does} \\ \textit{not always follow a loop.}
\end{center}

\begin{figure}[!htbp]
  \centering
	 \includegraphics[width=0.6\linewidth]{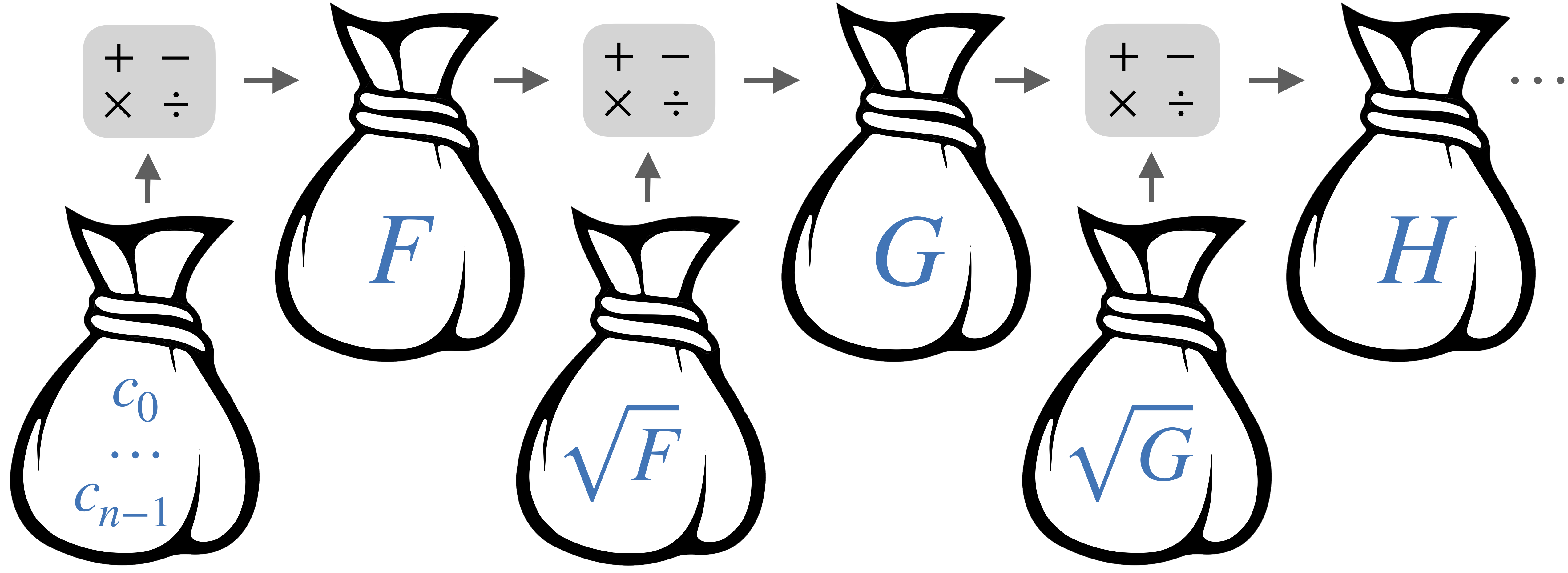}
   \caption{Ingredients used to build a formula. Combining coefficients $(c_0,...,c_{n-1})$ with $+,-,\times,\div$ defines an $F$-function. Combining $F$-functions and their roots $\sqrt{F}$ with $+,-,\times,\div$ defines a $G$-function, etc. \textit{See supplementary material for animated version.} \label{fig:FGH}}
\end{figure}

\textit{Formula ingredients.---}In this article we question the existence of a general \textit{formula} for the solutions of the general algebraic equation of degree $n$, equation \eqref{geneq}. By \textit{formula}, we mean some equality
\beq \label{formula}
s = \Phi(c_0,\ldots,c_{n-1}) \, ,
\eeq
where $s$ is a solution of equation \eqref{geneq} and $\Phi$ is some function of its coefficients $(c_0,\ldots,c_{n-1})$. The Abel--Ruffini theorem states that for $n\geq5$, no formula \textit{in terms of radicals} exists. ``In terms of radicals'' simply mean that the function $\Phi$ in equation \eqref{formula} can be constructed solely in terms of
\begin{itemize}
  \item the coefficients $(c_0,\ldots,c_{n-1})$ \,,
  \item the elementary operations $+$, $-$, $\times$, $\div$\ and $\psqrt$ \,.
\end{itemize}
Leaving $\psqrt$ aside, if we constrain ourselves to a formula combining the coefficients $(c_0,\ldots,c_n)$ and the four operations $+,-,\times,\div$, we obtain what we will call an $F$-formula, or simply an \textit{$F$-function}. Examples of such $F$-functions are
\beq
    F=1 \,, \quad F=-\frac{c_6}{2} \,,\quad F=c_8^2-7c_2 \,.
\eeq
They are the elementary building blocks for constructing formulae. In particular, they encompass integers, the coefficients themselves, as well as polynomials and rational functions of the coefficients. Clearly, if two coefficients each follow a loop simultaneously, then their sum, difference, product, and quotient also follow a loop. As they are built with only these four operations,  $F$-functions enjoy the same property. In other words:
\begin{center}
    \textit{When $(c_0,\ldots,c_{n-1})$ follow a loop,}\\ \textit{$F$-functions also follow a loop.}
\end{center}
This property of $F$-functions is not shared by $\sqrt{F}$-functions, i.e., expressions that are roots of $F$-functions, e.g., $\sqrt[7]{c_0}$ or $\sqrt[2]{1-3c_2}$. (recall the notation in the subsection ``\textit{Loops}''). In particular, if we denote by $G$-function a combination of $F$- and $\sqrt{F}$-functions together with $+,-,\times,\div$, then we have the following:
\begin{center}
    \textit{When $(c_0,\ldots,c_{n-1})$ follow a loop,}\\ \textit{$G$-functions do not always follow a loop.}
\end{center}
A $G$-function is a new type of ingredient as it may include expressions with one level of roots, such as
\beq
    G=-\frac{c_5}{2}+\frac{1}{2}\sqrt[2]{c_4^2-4c_1} \,.
\eeq
We can keep going like this to construct formulae with higher number of \textit{nested roots}, i.e., roots in roots. For example, we can combine $G$-functions and $\sqrt{G}$-functions with $+,-,\times,\div$ to make $H$-functions. These may contain up to two levels of nested roots, such as
\beq
  H=c_4-\sqrt[3]{7c_2}+\sqrt[5]{-\frac{c_0}{2}+\sqrt[5]{c_1^2-4c_6}} \,,
\eeq
and so on, as summarized in Figure \ref{fig:FGH}. With this nomenclature, we can make arbitrarily complex expressions involving $+,-,\times,\div$ and $\psqrt$, and at the same time keep track of the number of nested roots appearing in the formula. Conversely, any formula constructed with $+,-,\times,\div,\psqrt$ can be built using this procedure, provided that we look high enough in the ``$\ldots$'' of the list of ingredients $(F,G,H,\ldots)$.

We have now covered all the tools necessary: permutations of $(s_1,\ldots,s_n)$, loops, and $F,G,H$-functions. Let us now apply all these concepts to the degree $n$ equation, starting with $n=2$, to understand the Abel--Ruffini theorem when $n=5$.

\section{The quadratic equation}Our journey toward the Abel--Ruffini theorem starts with considerations of the much more familiar quadratic equation. In particular, considering only the case $n=2$, we will prove a first impossibility result, actually valid for $n\geq 2$. The ideas developed here are rather simple but also at the heart of the proof of the Abel--Ruffini theorem.

\textit{Vieta's formulae.---}Let us consider the general quadratic equation
\beq \label{equad}
z^2+c_1z+c_0=0 \, .
\eeq
As mentioned previously, the fundamental theorem of algebra informs us that this equation admits exactly two complex solutions $s_1$ and $s_2$. Let us then write it in the factored form $(z-s_1)(z-s_2) = 0$ and expand this product, ordering the terms by power of $z$. We find a new expansion that can be compared to equation \eqref{equad}. By identification, we obtain the so-called Vieta's formulae:
\beq \label{vietaquad}
    c_1 = - (s_1 + s_2) \quad \text{and} \quad c_0 = s_1 s_2 \, .
\eeq
This kind of relation between the coefficients and the solutions can be established for any degree $n\geq 2$. For example, equation \eqref{vietaquad} generalizes nicely to $c_{n-1}=-\sum_\ui s_\ui$ and $c_0=(-1)^n\prod_\ui s_\ui$, for any $n\geq 2$. In any case, these formulae always reveal the same fundamental property:

\begin{center}
    \textit{Coefficients $(c_0,\ldots,c_{n-1})$ are symmetric}\\ \textit{functions of the solutions $(s_1,\ldots,s_n)$.}
\end{center}

\begin{figure}[!htbp]
  \centering
	 \includegraphics[width=0.5\linewidth]{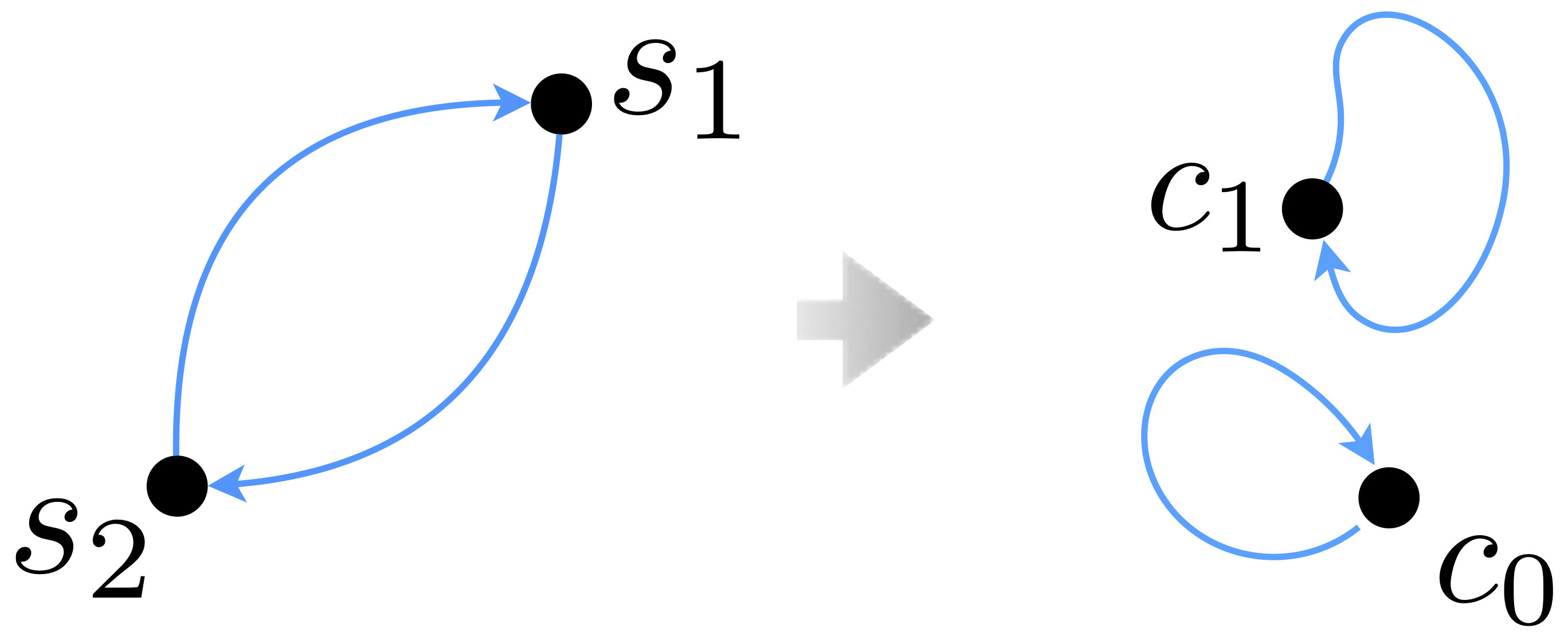}
   \caption{The transposition $(12)$ on the solutions $(s_1,s_2)$ induces a loop on the coefficients $(c_0,c_1)$. \textit{See supplementary material for animated version.} \label{fig:vieta}}
\end{figure}

In particular, for the $n=2$ case here at hand, if one permutes $s_1$ and $s_2$ by moving them continuously in the complex plane (using, for example, the transposition (12) depicted in Figure \ref{fig:permu}), then the coefficients $(c_0,c_1)$ will each move on some path, but eventually they must come back to their original location as they are symmetric in $(s_1,s_2)$. In other words, they will follow a \textit{loop}, as depicted in Figure \ref{fig:vieta}.

\begin{figure}[!htbp]
  \centering
	 \includegraphics[width=0.5\linewidth]{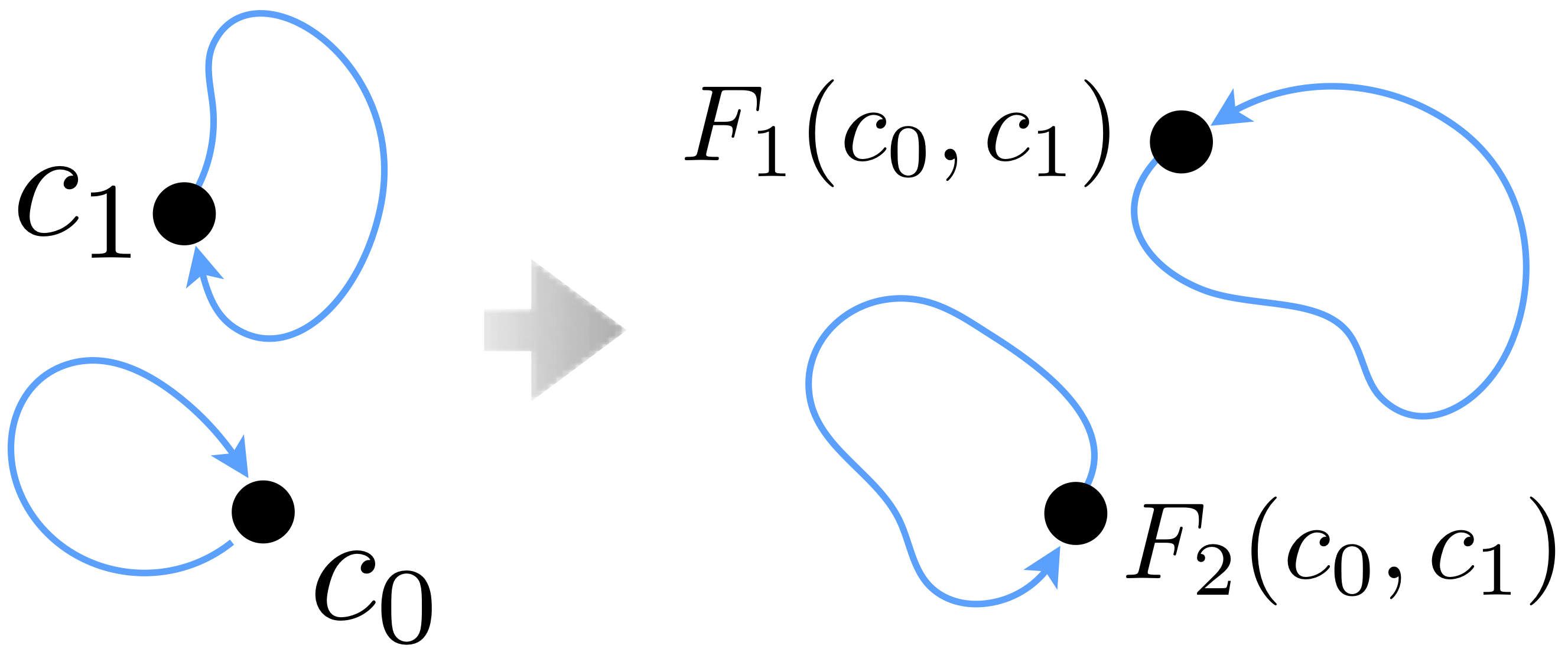}
   \caption{A loop followed by $(c_0,c_1)$ also induced a loop on the $F$-functions $F_1(c_0,c_1)$ and $F_2(c_0,c_1)$. \textit{See supplementary material for animated version.} \label{fig:permuquad}}
\end{figure}

\textit{A first impossibility result.---}Because it is the central idea at play, let us rephrase the symmetry in Vieta's formulae geometrically:
\begin{center}
    \textit{When $(s_1,\ldots,s_n)$ undergo a permutation}\\ \textit{$(c_0,\ldots,c_{n-1})$ each follow a loop.}
\end{center}
This remarkable fact can be used to obtain a first impossibility result, as follows. Suppose that the solutions $s_1$ and $s_2$ of the quadratic equation are given by two formulae of the type
\beq \label{quadformtest}
    s_1 = F_1(c_0,c_1) \quad \text{and} \quad s_2=F_2(c_0,c_1) \, ,
\eeq
with $F_1,F_2$ two $F$-functions (i.e., expressions involving $(c_0,c_1)$ and the symbols $+,-,\times,\div$). Now, picture $(s_1,s_2)$ and $(c_0,c_1)$ in the complex plane, and study the following process:
\begin{itemize}
    \item connect $s_1$ and $s_2$ with paths inducing the transposition $(12)$, and make them move along these paths (see Fig.~\ref{fig:vieta});
    \item as $s_1$ and $s_2$ move around, $c_0$ and $c_1$ each travel on a loop, as seen previously (see Fig.~\ref{fig:vieta}),
    \item while $c_0$ and $c_1$ follow their own loop, the two $F$-functions $F_1$ and $F_2$ will also follow a loop, as argued earlier (see Fig.~\ref{fig:permuquad}).
\end{itemize}

At the end of this process, $s_1$ and $s_2$ have permuted, yet both $F_1$ and $F_2$ have followed a loop. Consequently, $F_1$ and $F_2$ cannot equal $s_2$ and $s_2$ respectively, and no formula such as in \eqref{quadformtest} exists. This impossibility result actually holds for an equation of any degree $n\geq 2$. Indeed, it suffices to pick two of the $n$ solutions to the degree $n$ equation, name them $s_1$ and $s_2$, and apply the above recipe. The conclusion is thus:
\begin{center}
    \textit{Using only $F$-functions, no formula solving} \\ \textit{the general equation} \eqref{geneq} \textit{can be found for $n\geq 2$.}
\end{center}
This is our \textit{first impossibility result}. In particular, it means that we have no chance of finding a formula for the cubic equation using only $F$-functions either. To see which extra ingredients are needed, let us examine closely the well-known quadratic formula.

\textit{Discussion: the quadratic formula.---}The quadratic formula is derived most easily by ``completing the square'' in equation \eqref{equad} to get $(z+\tfrac{c_1}{2})^2=\tfrac{1}{4}c_1^2-c_0$.
Using our notation $\sqrt[2]{\phantom{\ui}}$ for any of the two square roots, we easily obtain a formula for the general solution $s$ of equation \eqref{equad} as
\beq \label{quadformula}
s = -\frac{c_1}{2} + \frac{1}{2}\sqrt[2]{c_1^2-4c_0} \, .
\eeq
This formula alone corresponds to two solutions, one for each of the two square roots on the right-hand side. Moreover, notice how this root indeed points to the same direction as our impossibility result: we need to add $\sqrt{F}$-function to the list of ingredients. One last note: just as the Abel--Ruffini theorem, the impossibility result just derived tells something about the \textit{general} quadratic equation. However, there exists \textit{some} quadratic equations with given, explicit coefficient that admit a formula in terms of $+,-,\times,\div$.

\section{The cubic equation}Let us now try to construct a formula for the solutions of the general cubic equation. The equation reads
\beq \label{cueq}
z^3+c_2z^2+c_1z+c_0=0 \, .
\eeq
Let $(s_1,s_2,s_3)$ be its three complex solutions. Learning from our previous findings, we now add $\sqrt{F}$-functions to the list of ingredients. Therefore, we assume that there exists some formulae of the type
\beq \label{rootcubic}
s_\ui = G_\ui(c_0,c_1,c_2) \, , \quad \text{for } \ui \in \{1,\ldots,3\} \, ,
\eeq
with $G_1,G_2,G_3$ three $G$-functions (combinations of $F$ and $\sqrt{F}$ with $+,-,\times,\div$). Our second impossibility result will consist in showing that such a formula cannot exist. Our previous method is not guaranteed to work: \textit{yes}, the coefficients still follow loops as solutions permute, but \textit{no}, $G$-functions do not generally follow loops in these circumstances, as we have already seen. We need to change our plan. \\

\begin{figure}[!htbp]
  \centering
	 \includegraphics[width=0.9\linewidth]{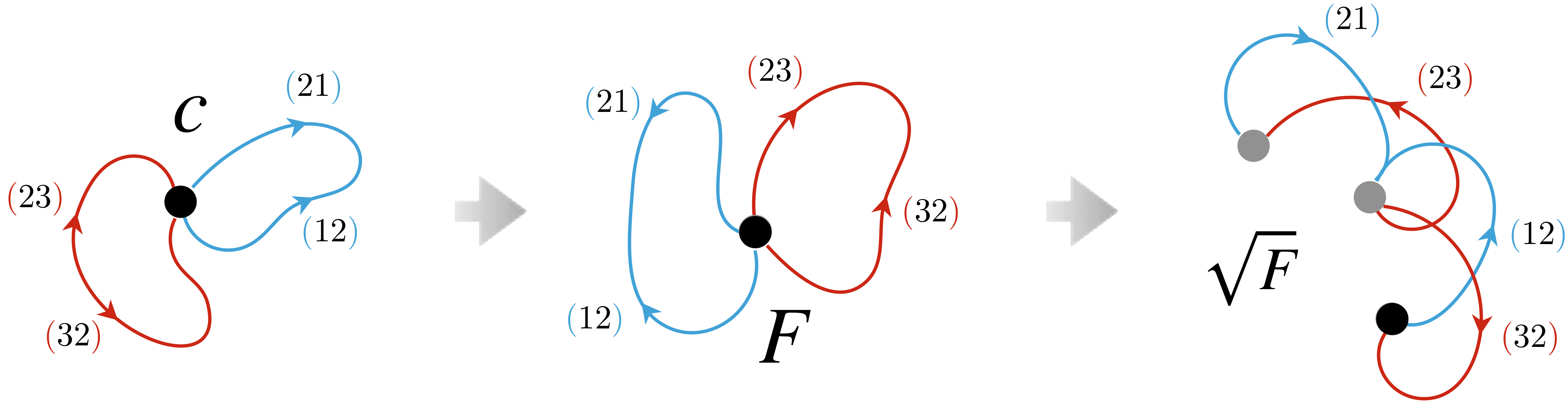}
   \caption{Effect of the commutator $[\cyan{(}12\cyan{)},\red{(}23\red{)}]$ on a coefficient $c$ (left), on an $F$-function (center) and on $\sqrt{F}$-function (right). After the process, both $F$ and $\sqrt{F}$ have followed a loop. Notice the loop followed by $\sqrt{F}$ consisting in four unclosed paths. \textit{See supplementary material for animated version.} \label{fig:permucubic}}
\end{figure}

\textit{Introducing commutators.---}Consider the transposition $\cyan{(}12\cyan{)}$ that induces a loop $\gamma_1$ on $F$ and thus an unclosed path $\omega_1$ on $\sqrt{F}$. Consider also $\red{(}23\red{)}$, inducing a loop $\gamma_2$ on $F$ and a path $\omega_2$ on $\sqrt{F}$. Now perform the following sequence of transpositions, called the \textit{commutator} of $\cyan{(}12\cyan{)}$ and $\red{(}23\red{)}$, and denoted
\beq  \label{permcubi}
[\cyan{(}12\cyan{)},\red{(}23\red{)}] = \cyan{(}12\cyan{)} \red{(}23\red{)} \cyan{(}12\cyan{)}^{-1} \red{(}23\red{)}^{-1} \, .
\eeq
Since $\cyan{(}12\cyan{)}^{-1}$ is simply $\cyan{(}21\cyan{)}$, and $\red{(}23\red{)}^{-1}=\red{(}32\red{)}$, it turns out that $[\cyan{(}12\cyan{)},\red{(}23\red{)}]$ is simply the cycle $(123)$. In fact, this is true with any pair of transposition, i.e.,
\beq \label{permu3}
[(ij),(jk)] = (ijk) \, .
\eeq
Therefore, $[\cyan{(}12\cyan{)},\red{(}23\red{)}]$ does permute the three solutions $(s_1,s_2,s_3)$. But what is its effect on numbers like $F$ and $\sqrt{F}$ ? Clearly, $F$ follows a sequence of loops $\gamma_1\gamma_2 \gamma_1^{-1}\gamma_2^{-1}$, which is itself a loop. The number $\sqrt{F}$, however, follows a sequence of unclosed paths $\omega_1 \omega_2\omega_1^{-1}\omega_2^{-1}$ (visiting other roots) but closes on itself by construction; see Figure \ref{fig:permucubic}.

\textit{Conclusion.---}With the permutation $(123)$ written as the commutator $[\cyan{(}12\cyan{)}, \red{(}23\red{)}]$, we reach the same conclusion as in the quadratic case: while $(s_1,s_2,s_3)$ undergoes the permutation $(123)$, both $F$ and $\sqrt{F}$ follow a loop (and thus any $G$-function). Consequently, there cannot be equalities given by \eqref{rootcubic}. Again, this holds for the general equation of degree $n\geq3$, too, as it suffices to pick up three solutions out of the $n\geq3$, label them $s_1,s_2,s_3$, and apply the above recipe. Therefore, we conclude:
\begin{center}
    \textit{Using only $G$-functions, no formula} \\ \textit{solving the general equation \eqref{geneq} can be found for $n\geq 3$.}
\end{center}
This is our \textit{second impossibility result}. We must emphasize that it works only if we apply the cycle $(123)$ as a commutator such as in equation \eqref{permcubi}. Had we just applied the cycle $(123)$ directly (i.e., without writing it as a commutator), there would have been no guarantee that $\sqrt{F}$ followed a loop. It is the commutator that allows us to \textit{discard} one level of roots, and thus $\sqrt{F}$, from the list of ingredients. Let us now put this new impossibility result to the test, by solving explicitely the cubic equation. \\

\textit{Discussion: the cubic formula.---}We follow the classical method found by Italian mathematicians of the sixteenth century. First, perform the change of variables $Z=z+c_2/3$, which ``removes'' the $z^2$ term in equation \eqref{cueq}, transforming it into
\beq \label{redcubic}
Z^3+3PZ+2Q=0 \, ,
\eeq
where $P=\tfrac{c_1}{3}-\tfrac{c^2_2}{9}$ and $Q=\tfrac{c_0}{2}+\tfrac{c_2^3}{27}-\tfrac{c_1c_2}{6}$.
Notice that both $P$ and $Q$ are $F$-functions of $(c_0,c_1,c_2)$. To solve equation \eqref{redcubic}, one then writes $Z=v+w$, where $v,w$ are two complex numbers to be chosen freely later on. Then, equation \eqref{redcubic} becomes $v^3+w^3+3(vw+P)(v+w)+2Q=0$, from which we can remove the second term by imposing that $v,w$ satisfy $vw=-P$. By cubing the latter, we then obtain two equations for two unknowns, namely
\beq
	v^3+w^3=-2Q \quad \text{ and } \quad v^3w^3=-P^3 \, .
\eeq
These equations can be solved simultaneously for $v^3$ and $w^3$, since they explicitly give their sum and product, respectively. (These are nothing but Vieta's formulae for $n=2$; See equation \eqref{vietaquad}.) Using the quadratic formula, one obtains $v^3$ and $w^3$ in terms of $P$ and $Q$, takes their cube root and adds the result to obtain $v+w=Z$. Going back to the original unknown $z=Z-c_2/3$ gives the famous ``cubic formula''
\beq \label{cubicformula}
s = -\frac{c_2}{3} + \sqrt[3]{-Q + \sqrt{Q^2+P^3}} + \sqrt[3]{-Q - \sqrt{Q^2+P^3}}  \, .
\eeq
This formula gives three solutions $(s_1,s_2,s_3)$, one for each cube root. It is clear that this expression involves more than $F$ and $\sqrt{F}$ functions: indeed, the two cube roots are actually $\sqrt{G}$-functions. In a sense, the cubic formula above contains ``two levels'' of roots, whereas $G$-functions can only contain one, by definition. This kind of expression is thus called a \textit{nested root}. Our ``commutator trick'' was only able to remove \textit{one} level of roots. Perhaps two levels of commutators will remove two? If so, then it looks like a pattern is emerging\ldots

\section{The quartic equation}We now turn to the quartic equation
\beq \label{quartic}
z^4+c_3z^3+c_2z^2+c_1z+c_0=0 \, .
\eeq
For the cubic, we saw that $G$-functions are not enough to construct a formula, as we also needed $\sqrt{G}$ functions. Therefore, we start by assuming the existence of some formula for the four solutions
\beq \label{rootquartic}
s_\ui = H_\ui(c_0,c_1,c_2,c_3) \, , \quad \text{for } i\in\{1,\ldots,4\} \, .
\eeq
As before, the four functions $H_\ui$ are assumed to be $H$-functions, i.e., $G$- and $\sqrt{G}$-functions combined with $+,-,\times,\div$. As should be clear by now, it turns out that even with the extra ingredient $\sqrt{G}$, no general quartic formula can be constructed. Once again we will prove this by constructing an appropriate permutation of $(s_1,s_2,s_3,s_4)$. \\

\textit{A brief checkpoint.---}Once again, just as the first method did not work for cubic equations, the method used for the cubic case is not guaranteed to work for quartic equations either. Indeed, the commutator of transpositions induced a loop on $F$ and $\sqrt{F}$ (and thus on $G$). But a loop on $G$ generally does not induce a loop on $\sqrt{G}$, as we have seen many times. A summary of these previous methods is given on Table \ref{table:summary}.

\begin{table}[!htbp]
	\centering
			\begin{tabular}{|c|c|c|}
			\hline
			ingredient &\, $F$-functions \, & \, $G$-functions \, \\ \hline
			nested roots &\, $0$ \, & \, $1$ \, \\ \hline
			discarded by & \, \cyan{transpositions}	& \, \red{commutator} of \cyan{transpositions} \, \\ \hline
			with the path & \, \cyan{(}12\cyan{)} \, & \, \red{[}\cyan{(}12\cyan{)},\cyan{(}23\cyan{)}\red{]} = (123) \, \\ \hline
      for degree & \, $n \geq 2$ \, & \, $n \geq 3$ \, \\ \hline
			\end{tabular}
		\caption{Summary of the methods used to prove the first two impossibility results. \label{table:summary}}

\end{table}

But now a natural solution presents itself: what if we take the commutator of, say, $(123)$ and $(234)$, written as commutators themselves, using equation \eqref{permu3}? Let us examine this in detail. \\

\textit{Commutators, yet again.---}First we need to check that the commutator of $(123)$ and $(234)$ does indeed permute the four solution $(s_1,s_2,s_3,s_4)$. Fortunately it does, as a quick check reveals that
\beq \label{permu4}
[\cyan{(}123\cyan{)}, \red{(}234\red{)}]=(14)(23) \, ,
\eeq
which is a particular case of the more general formula $[(ijk),(jk\ell)] = (i\ell)(jk)$. Therefore, our commutator $[\cyan{(}123\cyan{)}, \red{(}234\red{)}]$ does indeed permute $(s_1,s_2,s_3,s_4)$. Now, let us examine how it affects $G$- and $\sqrt{G}$-functions, one step at a time:

\textbullet \,\, First, we apply the cycles $\cyan{(}123\cyan{)}=\cyan{[}(12),(23)\cyan{]}$ then $\red{(}234\red{)}=\red{[}(23),(34)\red{]}$. Since they are commutators, $G$-functions will follow two loops $\gamma_1,\gamma_2$ successively, coming back to their original positions. However, quantities like $\sqrt{G}$ will move along two (generally unclosed) paths $\omega_1$ and $\omega_2$. All this is exactly as in the cubic case. \\

\textbullet \,\, Second, we apply these two paths backwards, in reverse i.e., $\red{(}432\red{)}=\red{[}(43),(32)\red{]}$ and then $\cyan{(}321\cyan{)}=\cyan{[}(32),(21)\cyan{]}$. During these two, $G$-functions will follow $\gamma_2^{-1}\gamma_1^{-1}$, i.e. the previous loops backwards. Similarly, $\sqrt{G}$-functions will travel along $\omega_2^{-1}\omega^{-1}_1$.\\

What just happened is exactly the same as in the cubic case, except with $G$-functions in place of $F$-functions. In particular, $G$-functions follow the loop $\gamma=\gamma_1 \gamma_2 \gamma_1^{-1} \gamma_2^{-1}$; and $\sqrt{G}$-functions a sequence of unclosed paths $\omega_1 \omega_2 \omega_1^{-1} \omega_2^{-1}$, which closes on itself by construction. In other words, both $G$- and $\sqrt{G}$-functions followed a loop and thus any $H$-function will, too.
Our conclusion has therefore been reached: while $(s_1,s_2,s_3,s_4)$ undergoes the permutation $(14)(23)$ written as a commutator of commutators, any $H$-function follows a loop. Consequently, no formula \eqref{rootquartic} can exist. This result extends to any equation of degree $n\geq4$, as before, and constitutes our \textit{third impossibility result:}
\begin{center}
    \textit{Using only $H$-functions, no formula solving} \\ \textit{the general equation \eqref{geneq} can be found for $n\geq 4$.}
\end{center}
In particular, we can extend Table \ref{table:summary} with an additional column for the new ingredient, $H$-functions.
\begin{table}[!htbp]
	\centering
			\begin{tabular}{|c|c|}
			\hline
			ingredient &        \, $H$-functions \, \\ \hline
			nested roots &        \, 2 \, \\ \hline
			discarded by &      \, \green{commutator} of \red{commutator} of \cyan{transpositions}	\, \\ \hline
			with the path &          \, \green{[}\red{[}\cyan{(}12\cyan{)},\cyan{(}23\cyan{)}\red{]},\red{[}\cyan{(}23\cyan{)},\cyan{(}34\cyan{)}\red{]}\green{]} = (14)(23) \, \\ \hline
            for degree &        \, $n \geq 4$ \, \\ \hline
			\end{tabular}
		\caption{Extension of Table \ref{table:summary} to $H$-functions, for degree $n\geq4$. \label{table:summary4}}
\end{table}

\textit{Discussion: the quartic formula.---}As for the cubic case, our impossibility result does not imply that there is no quartic (nor quintic) formula. It just means that to construct one, one needs at least three levels of nested roots, and $H$-functions contain only two. It turns out that the quartic equation can be solved as follows and, indeed, it involves three levels of nested roots. As for the cubic case, we start by removing the $z^3$ term by the change of variables $Z=z+c_3/4$.
This brings equation \eqref{quartic} into the form
\beq \label{quart}
Z^4+PZ^2+QZ+R = 0 \, ,
\eeq
where $P,Q,R$ are three $F$-functions of $(c_0,c_1,c_2,c_3)$, whose expressions are long, but easily obtained. The next step is to transform equation \eqref{quart} into one that is quadratic in $Z^2$. For now, nothing guarantees that $PZ^2+QZ+R$ is a perfect square, but if it were, then equation \eqref{quart} could be factored into two equations quadratic in $Z^2$. One way is to write $Z^4$ in the equivalent form $Z^4=(Z^2+A)^2-2AZ^2-A^2$, for some complex number $A$ to be chosen freely later on. Inserting this in equation \eqref{quart} gives
\beq \label{interA}
(Z^2+A)^2+(P-2A)Z^2+QZ+R-A^2 = 0 \, .
\eeq
Now we can choose $A$ in equation \eqref{interA} such that the quadratic part $(P-2A)Z^2+QZ+R-A^2$ has the form of a perfect square. This will be the case if its discriminant $Q^2+4(P-2A)A^2$ vanishes. The latter amounts to
\beq \label{cubicA}
8A^3 -4PA^2  -Q^2  = 0 \, ,
\eeq
which is a cubic equation $A$. It can be solved using the cubic formula, giving a value of $A$ in terms of $P$ and $Q$ that is an $H$-function (recall the cubic formula \eqref{cubicformula}). Once $A$ takes this special value, equation \eqref{interA} becomes $(Z^2+A)^2+(P-2A)(Z-A)^2=0$, which can be factored easily into two quadratic polynomials in $Z$. The latter equations are solved easily using the quadratic formula. Since $A$ is an $H$-function, the solution for $Z$ will necessarily involve some $\sqrt{H}$ quantities, something which we did not include in equation \eqref{rootquartic}. This confirms our impossibility result, once again.
\section{The quintic equation}It seems at this point that things are becoming repetitive, and that a clear pattern emerges. For $n=2,3,4$, commutators could be used to reject formulae with too few nested roots in their expressions. However, we were still be able to solve the equation simply by allowing more levels of roots. But at $n=5$, this all breaks down, and this is why the quintic equation is a {very} special case. The goal of this section is to apply our methods to the case of degree $n=5$ and understand why it allows, not only to discard $4$ levels of nested roots (i.e., one more that the quartic case), but actually \textit{any} number of roots. \\

Let us pretend that we found a quintic formula, e.g.,
\beq \label{quintform}
s_\ui = \Phi_\ui(c_0,\ldots,c_4), \quad \text{for } i\in\{1,\ldots,5\} \,,
\eeq
with the five functions $\Phi_\ui$ built out of $H$-- and $\sqrt{H}$--functions.
If we follow the previous methods, summarized in Tables \ref{table:summary} and \ref{table:summary4}, it should be clear that (1) all $H$-functions will follow a loop from a commutator of commutators of the solutions (as in the quartic case), but (2) we will need one more level of commutators for the $\sqrt{H}$ terms.

As we now have five solutions to play with, let us consider for example the permutations $(123)$ and $(345)$ to construct a first commutator $[(123),(345)]$. An easy check shows that the latter is equal to $(235)$, and this commutator therefore permutes three of our solutions. In general, the following result holds at $n=5$:
\beq \label{3permu}
[(ijk),(k\ell m)] = (jkm) \, .
\eeq
But now, contrary to the previous cases, we have something rather remarkable with equation \eqref{3permu}. It shows that any cycle $(jkm)$ can be written as a commutator of two other cycles, namely $[(ijk),(k\ell m)]$. But notice that this is true for \textit{any} cycle $(jkm)$, including $(ijk)$ and $(k\ell m)$ on the left-hand side of equation \eqref{3permu} itself. In other words, this formula can be applied to itself, again and again, allowing us to write $(jkm)$ as a commutator of as many commutators as needed. Since a number $N\in\NN$ of commutators allows us to discard precisely $N$ levels of roots in a formula (see Tables \ref{table:summary} and \ref{table:summary4}), we can actually discard \textit{any} number of roots in any candidate
quintic formula. The Abel--Ruffini theorem follows immediately from this remark, but let us give a more detailed explanation.

Suppose that, in the quintic formula, equation \eqref{quintform}, we use a $\Phi$-function made of $+,-,\times,\div$, along with $N$ levels of roots, for some $N\in\NN$. To construct this $\Phi$-function, we have at our disposal several ingredients: $F$-functions (no roots), $G$-functions (one level of root), $H$-functions (two levels of root), and so on. As always, we start by choosing a permutation of the solutions, say $(123)$, that discards any $F$-functions (no roots). Next, using equation \eqref{3permu}, we write $(123)$ as a commutator, for example:
\beq \label{3permu2}
\cyan{(} 123 \cyan{)}=\cyan{[} \red{(}412\red{)},\orange{(}253\orange{)} \cyan{]}\, .
\eeq
When applied to $(s_1,\ldots,s_5)$, this commutator discards the $G$-functions from the list of ingredients (one level of roots). Now we keep going: we write the cycles $(412),(253)$ appearing in equation \eqref{3permu2} as commutators themselves, again using equation \eqref{3permu}. We obtain $(123)$ expressed with two commutators:
\beq \label{123dooublecommu}
\cyan{(} 123 \cyan{)}=\cyan{[} \red{[}(341),(152)\red{]} , \orange{[}(425),(513)\orange{]} \cyan{]} \, ,
\eeq
which removes $H$-functions (expression with two levels of roots). By iterating equation \eqref{3permu} $N-2$ more times, we end up writing $(123)$ as a combination of $N$ commutators. When the latter is applied to $(s_1,\ldots,s_5)$, the solutions $s_1,s_2,s_3$ will permute; and yet any expression of the coefficients with $N$ or less roots will follow a loop. Since a $\Phi$-function is made up of all these ingredients, $\Phi_1,\Phi_2,\Phi_3$ go back to their original position. Clearly this contradicts equation \eqref{quintform}. This result generalizes to an equation of any degree $n\geq 5$ by picking five of its solutions, as before. Moreover, since $N$ is arbitrary, we conclude that no number of roots will be sufficient to write a formula. Our conclusion is therefore:
\begin{center}
\textit{No formula exists for the solution to the general} \\ \textit{equation of degree five or more, using only} \\ \textit{the operations $+,-,\times,\div,$ and $\psqrt$} \,,
\end{center}
i.e., the Abel--Ruffini theorem itself. A last remark is in order. Why the fifth degree, and not the fourth or sixth? This all boils down to the possibility of writing a permutation of at least two solutions as a sequence of commutators. A formula such as equation \eqref{3permu} can only be iterated indefinitely when it involves five or more elements. For four or fewer elements, any sequence of commutators of transposition and/or cycles will necessarily end, i.e., end up being the trivial permutation that ``does nothing.'' The reader familiar with group theory will here recognize the notions of perfect or solvable group.
\section{Conclusions}To conclude this article, we would like to first make some comments on the various advantages and disadvantages of this proof, compared to the usual proof using Galois theory. First of all, the present proof does not say that no equations of degree five or higher can be solved; but only that a general formula (valid for the general equation) cannot be written using only $+,-,\times,\div,$ and $\psqrt$. Indeed, some equations of degree $n\geq 5$ can actually be solved explicitly (see \cite{SpWi} for a nice and thorough exposition on the quintics $z^5+az+b=0$ that are solvable by radicals.) Galois theory, on the other hand, is perfectly able to say whether a given equation is solvable or not. On the other hand, the present proof can be extended to also account for continuous (and single-valued) functions of the coefficients (such as $\exp, \sin,\ldots$) in the list of ingredients. Indeed, just like $+,-,\times,\div$, these functions follow a loop when the coefficients do. Galois theory is unable to provide for this, as it only accounts for expressions in terms of radicals.\\

We hope that the present proof will be seen not only as a simplified and elementary demonstration of the Abel--Ruffini's theorem, but also as a complementary result, as it helps to explain why the $n=5$ case is so special, and why the quadratic, cubic, and quartic formulae have such ``nested roots'' structures. It is also a good and instructive exercise to complete the present proof sketch with rigorous arguments (we encourage students to give it a try!). For more insight on this topic, one should definitely watch Katz's video \cite{Katz} and read Goldmakher's paper \cite{Goldmakher}, which have both inspired this article.
We end this paper by providing more references that should help the interested reader to get started with topics that are based on (and broadly extend) the ideas presented here: (1) an interactive blog article by F.~Akalin \cite{Akalin}; (2) an article by H.~Zoladek \cite{Zoladek} that deals with similar but more advanced ideas; (3) the original book of Alekseev mentioned in the introduction \cite{Alekseev} on the Abel--Ruffini Theorem (see \url{www.maths.ed.ac.uk/~v1ranick/papers/abel.pdf} for a free digital copy). Devised as a problems-and-solutions book, it discusses many advanced concepts in a very pedagogical and extremely well-written manner.

\end{document}